\documentclass[final,3p, times]{elsarticle}

\usepackage{amssymb,amsmath,color}
\usepackage[mathscr]{eucal}
 \usepackage{lineno}

 \journal{ }
\renewcommand{\journal}{  TBD }
\newtheorem{thm}{Theorem}[section]
\newtheorem{prop}[thm]{Proposition}

\newproof{pf}{Proof}

\begin{document}

\begin{frontmatter}



\title{Metric dimension of some distance-regular graphs}

\author[JG]{Jun Guo}\ead{guojun$_-$lf@163.com}
\author[KW]{Kaishun Wang\corref{cor}}
\ead{wangks@bnu.edu.cn}
\author[FL]{Fenggao Li}\ead{fenggaol@163.com}
\cortext[cor]{Corresponding author}

\address[JG]{Math. and Inf. College, Langfang Teachers' College, Langfang  065000,  China }
\address[KW]{Sch. Math. Sci. \& Lab. Math. Com. Sys.,
Beijing Normal University, Beijing  100875, China}
\address[FL]{Dept. of Math., Hunan Institute of Science and Technology, Yueyang  414006,   China}

\begin{abstract}
A resolving set of a graph is a set of vertices
with the property that the list of distances from any vertex to those in the
set uniquely identifies that vertex. In this paper, we construct
a resolving set of Johnson graphs, doubled Odd graphs, doubled Grassmann graphs
and twisted  Grassmann graphs, respectively, and obtain
the upper bounds on the metric dimension of these graphs.
\end{abstract}

\end{frontmatter}

\section{Introduction}

Let $\Gamma$ be a connected graph. For any two vertices $u$ and $v$,
$d(u, v)$ denotes the distance between $u$ and $v$.
A {\it resolving set} of $\Gamma$ is a set of vertices $S = \{v_1,\ldots,v_k\}$
such that, for each vertex $w$, the ordered list of distances $D(w|S)=(d(w, v_1),\ldots,d(w, v_k))$
uniquely determines $w$. That is, $S$ is a resolving set of $\Gamma$
if, for any two distinct vertices $u$ and $v$, $D(u|S)\not= D(v|S)$.
The {\it metric dimension} of $\Gamma$, denoted by $\mu(\Gamma)$,
is the smallest size of all resolving sets of $\Gamma$.

Metric dimension is a well-known parameter in graph theory. It was
first introduced in the 1970s, independently by Harary and Melter
\cite{HM} and by Slater \cite{Sla}. Computing the metric dimension
of a graph is an NP-hard problem \cite{BC}. In recent years, a
considerable literature has developed in graph theory.
 C${\rm\acute{a}}$ceres et al. \cite{CHMPPSW} studied the metric
dimension of Cartesian products of graphs.  Chartrand et al.
\cite{CEJO} determined all connected graphs of order $n$ having
metric dimension 1, $n-2$ or $n-1$, and presented a new proof on the
metric dimension of a tree. An interesting case is that of
distance-regular graphs. Bailey et al. \cite{BCGG} obtained an upper
bound on the metric dimension of Johnson graphs. Bailey and Meagher
\cite{BM} constructed a resolving set of Grassmann graphs.
Chv${\rm\acute{a}}$tal \cite{Chv} obtained an upper bound on the
metric dimension of Hamming graphs. Feng and Wang  \cite{FW}
obtained an upper bound on the metric dimension of bilinear forms
graphs.

In this paper, we construct
a resolving set of Johnson graphs, doubled Odd graphs, doubled Grassmann graphs
and twisted  Grassmann graphs, respectively, and obtain
the upper bounds on the metric dimension of these graphs.

\section{Johnson graphs}
For any positive integer $m$, let $[m]=\{1,2,\ldots,m\}$  and let
$\binom{[m]}{i}$ be the set of all $i$-subsets of $[m]$. Given
integers $2\leq 2e\leq n,$ the Johnson graph $J(n,e)$ has the vertex
set $\binom{[n]}{e}$ such that two vertices $P$ and $Q$ are adjacent
if and only if $|P\cap Q|=e-1$. Then $J(n,e)$ is a distance-regular
graph with diameter $e$ (see \cite{BCN}). For any two vertices $P$
and $Q$, we have $d(P,Q)=i$ if and only if $|P\cap Q|=e-i$. For
$e=1$, $J(n,e)$ is a complete graph. It is obvious that
$\mu(J(n,1))= n-1$. For $e=2$, by \cite[Corollary~3.33]{BC},
$\mu(J(3,2))=2,\mu(J(4,2))=3,\mu(J(5,2))=3$ and
$\mu(J(n,2))=\frac{2}{3}(n-i)+ i$, where $n\geq6$ and $n\equiv i$
(mod 3), $i\in\{0, 1, 2\}$. For $e\geq3$, Bailey et al. obtained the
following result.

\begin{prop}{\rm(\cite{BCGG})}\label{prop3.1}
Let $e\geq3$. Then $\mu(J(n,e))\leq (e+1)\lceil n/(e+1)\rceil.$
\end{prop}

For $n=2e+1$, we improve the bound of Proposition~\ref{prop3.1} and
obtain the following result.

\begin{thm}\label{thm1.1.j}
Let $e\geq3$. Then $\mu(J(2e+1,e))\leq 2e$.
\end{thm}

\begin{pf}
Let $X=[e]$ and $Y=[2e]\backslash X$.
To prove
$${\cal M}=\left\{M_1\cup\{2e+1\}\mid M_1\in\binom{X}{e-1}\right\}\cup\left\{M_2\cup \{2e+1\}\mid M_2\in\binom{Y}{e-1}\right\}$$
is a resolving set, it is sufficient to
show that, for any two distinct vertices $P$ and $Q$,
there exists a vertex $W\in{\cal M}$ such that
$|P\cap W|\not=|Q\cap W|$.

{\it Case} 1: $P\cap[2e]\not=Q\cap[2e]$. Then there exists
a $Z\in\{X,Y\}$ such that $P\cap Z\not= Q\cap Z$. Suppose
$s=|P\cap Z|\leq|Q\cap Z|=r$. Then $1\leq r\leq e$.

{\it Case} 1.1: $r=e$. Then $Q=Z$. Since $P\cap Z\not= Q\cap Z$,
$s<e$. If $s=0$, take an $(e-1)$-subset $W_1$ of $Z$. Then
$W=W_1\cup\{2e+1\}\in{\cal M}$ and $|P\cap W|\leq1<e-1=|Q\cap W|$.
If $1\leq s\leq e-2$, let $W_1$ be an $(e-1)$-subset of $Z$ satisfying $|P\cap W_1|=s-1$.
Then $W=W_1\cup\{2e+1\}\in{\cal M}$ and $|P\cap W|\leq s< e-1=|Q\cap W|.$
If $s=e-1$ and $2e+1\in P$, then $W=(P\cap Z)\cup\{2e+1\}\in{\cal M}$ and
$|P\cap W|= e> e-1=|Q\cap W|.$ If $s=e-1$ and $2e+1\not\in P$,
let $W_1$ be an $(e-1)$-subset of $Z$ satisfying $|P\cap W_1|=s-1$.
Then $W=W_1\cup\{2e+1\}\in{\cal M}$ and $|P\cap W|= s-1< e-1=|Q\cap W|.$

{\it Case} 1.2: $r=1$. If $s=0$ and either $2e+1\in Q$ or $2e+1\not\in P\cup Q$,
let $W_1$ be an $(e-1)$-subset of $Z$
containing $Q\cap Z$. Then $W=W_1\cup\{2e+1\}\in{\cal M}$ and $|P\cap W|<|Q\cap W|$.
If $s=0,2e+1\in P$ and $2e+1\not\in Q$,
take $W_1=Z\backslash(Q\cap Z)$.
Then $W=W_1\cup\{2e+1\}\in{\cal M}$ and $|P\cap W|=1>0=|Q\cap W|$.
If $s=1$ and $2e+1\in P$, let $W_1=Z\backslash(Q\cap Z)$.
Then $W=W_1\cup\{2e+1\}\in{\cal M}$ and $|P\cap W|=2>1\geq|Q\cap W|$.
If $s=1$ and $2e+1\not\in P$, let $W_1=Z\backslash(P\cap Z)$.
Then $W=W_1\cup\{2e+1\}\in{\cal M}$ and $|P\cap W|=0<1\leq|Q\cap W|$.

{\it Case} 1.3: $1<r<e$. If $0\leq s\leq r-2$, take an $(e-1)$-subset $W_1$
of $Z$ containing $Q\cap Z$. Then
$W=W_1\cup\{2e+1\}\in{\cal M}$ and $|P\cap W|\leq s+1<r\leq|Q\cap W|$.
If $s=r-1$ and either $2e+1\in Q$ or $2e+1\not\in P\cup Q$,
let $W_1$ be an $(e-1)$-subset of $Z$ containing $Q\cap Z$.
Then $W=W_1\cup\{2e+1\}\in{\cal M}$ and $|P\cap W|<|Q\cap W|.$
If $s=r-1,2e+1\in P$ and $2e+1\not\in Q$,
let $W_1$ be an $(e-1)$-subset of $Z$ satisfying $P\cap Z\subseteq W_1$
and $|W_1\cap Q\cap Z|=r-1$.
Then $W=W_1\cup\{2e+1\}\in{\cal M}$ and $|P\cap W|= s+1>r-1=|Q\cap W|.$
If $s=r$ and $2e+1\in Q$,
let $Q\cap Z=\{\alpha_1,\ldots,\alpha_r\}$.
Since $P\cap Z\not=Q\cap Z$,
there exists a $\beta\in (P\cap Z)\backslash (Q\cap Z)$ such that
$\{\alpha_1,\ldots,\alpha_r,\beta\}$ is an $(r+1)$-subset.
Let $Z=\{\alpha_1,\ldots,\alpha_r,\beta,\gamma_1,\ldots,\gamma_{e-r-1}\}$
and $W_1=\{\alpha_1,\ldots,\alpha_r,\gamma_1,\ldots,\gamma_{e-r-1}\}$.
Then $W=W_1\cup\{2e+1\}\in{\cal M}$. By $\beta\not\in W$ and $\beta\in P\cap Z$,
$|P\cap W|\leq|P\cap Z|=r<r+1=|Q\cap W|.$
If $s=r$ and $2e+1\in P$, then there exists a $W\in{\cal M}$ such that
$|P\cap W|=r+1>r\geq|Q\cap W|.$
If $s=r$ and $2e+1\not\in P\cup Q$, then there exists an $(e-1)$-subset
$W_1$ such that $|P\cap W_1|<|Q\cap W_1|$. It follows that $W=W_1\cup\{2e+1\}\in{\cal M}$
and $|P\cap W|<|Q\cap W|.$

{\it Case} 2: $P\cap[2e]=Q\cap[2e]$. Then $P\cap[e]=Q\cap[e]$.
Let $W_1$ be an $(e-1)$-subset of $[e]$. Then
$W=W_1\cup\{2e+1\}\in{\cal M}$ and
$|P\cap W|\not=|Q\cap W|$.

Hence, ${\cal M}$ is a resolving set. Finally, we obtain the bound by observing that
$|{\cal M}|=2e.$
\qed \end{pf}

\section{Doubled Odd graphs and their $q$-analogue}

The doubled Odd graph, denoted by $O(2e+1,e,e+1)$,
is a bipartite graph with bipartition $\binom{[2e+1]}{e}\cup\binom{[2e+1]}{e+1}$ such that
two vertices $P\in\binom{[2e+1]}{e}$ and $Q\in\binom{[2e+1]}{e+1}$
are adjacent if and only if $P\subseteq Q$.
Then $O(2e+1,e,e+1)$ is a distance-regular graph
with diameter $2e+1$ (see \cite{BCN}).
For any two vertices $P$ and $Q$, we have
$d(P,Q)=2i+|\,|P|-|Q|\,|$ if and only if
$|P\cap Q|=\min\{|P|,|Q|\}-i$. For $e=1$, $O(3,1,2)$ is a hexagon.
It is obvious that $\mu(O(3,1,2))= 2e$.
For $e\geq2$, we obtain the following result.

\begin{thm}\label{thm1.0.1}
Let $e\geq2$. Then $\mu(O(2e+1,e,e+1))\leq 2e+1$.
\end{thm}

\begin{pf}
We give an explicit construction of a resolving set.

Let $Y=[2e+1]\backslash[e+1]$,
${\cal W}=\{W\subseteq [e+1]\mid |W|=e\}$
and ${\cal Y}=\left\{W\cup\{1\} \mid W\in \binom{Y}{e-1}\right\}$.
We will show that ${\cal M}={\cal W}\cup{\cal Y}$ is a resolving set.
Since $O(2e+1,e,e+1)$ is bipartite, we only need to show that,
for any two distinct vertices $P$ and $Q$ in the same part,
there exists a vertex $W\in{\cal M}$ such that
$|P\cap W|\not=|Q\cap W|$.

Suppose $P,Q\in\binom{[2e+1]}{e+i}$, where $i=0$ or 1.

{\it Case} 1: $P\cap [e+1]\not=Q\cap [e+1]$. Suppose
$s+i=|P\cap [e+1]|\leq|Q\cap [e+1]|=r+i$. Then $1-i\leq r\leq e$.

{\it Case} 1.1: $s<r$. For $i=0$, let $W$ be an $e$-subset of $[e+1]$ containing $Q\cap[e+1]$.
Then $W\in{\cal M}$ and $$|P\cap W|\leq|P\cap [e+1]|=s<r=|Q\cap W|.$$
For $i=1$, let $W$ be an $e$-subset of $[e+1]$ satisfying $|P\cap W|=s$.
Then $W\in{\cal M}$ and $|P\cap W|=s<r\leq|Q\cap W|$.

{\it Case} 1.2: $s=r$. Let $Q\cap [e+1]=\{\alpha_1,\ldots,\alpha_r,\alpha_{r+i}\}$.
Since $P\cap [e+1]\not=Q\cap [e+1]$,
there exists a $\beta\in (P\cap [e+1])\backslash (Q\cap [e+1])$ such that
$\{\alpha_1,\ldots,\alpha_r,\alpha_{r+i},\beta\}$ is an $(r+1+i)$-subset of $[e+1]$.
Let $[e+1]=\{\alpha_1,\ldots,\alpha_r,\alpha_{r+i},\beta,\gamma_1,\ldots,\gamma_{e-r-i}\}$
and $W=\{\alpha_1,\ldots,\alpha_r,\alpha_{r+i},\gamma_1,\ldots,\gamma_{e-r-i}\}$.
Then $W\in{\cal M}$. By $\beta\not\in W$ and $\beta\in P\cap [e+1]$,
$$|P\cap W|<|P\cap [e+1]|=|Q\cap [e+1]|=r+i=|Q\cap W|.$$

{\it Case} 2: $P\cap [e+1]=Q\cap [e+1]$. Then $P\cap Y\not=Q\cap Y$
and $|P\cap Y|=|Q\cap Y|$.
Suppose $|P\cap Y|=|Q\cap Y|=r$. Then $0<r< e$.
Let $Q\cap Y=\{\alpha_1,\ldots,\alpha_r\}$.
Since $P\cap Y\not=Q\cap Y$,
there exists a $\beta\in (P\cap Y)\backslash (Q\cap Y)$ such that
$\{\alpha_1,\ldots,\alpha_r,\beta\}$ is an $(r+1)$-subset of $Y$.
Let $Y=\{\alpha_1,\ldots,\alpha_r,\beta,\gamma_1,\ldots,\gamma_{e-r-1}\}$
and $W_1=\{\alpha_1,\ldots,\alpha_r,\gamma_1,\ldots,\gamma_{e-r-1}\}$.
By $\beta\not\in W_1$ and $\beta\in P\cap Y$,
$$|P\cap W_1|<|P\cap Y|=|Q\cap Y|=r=|Q\cap W_1|.$$
Therefore, $W=W_1\cup\{1\}\in{\cal M}$ and
 $|P\cap W|\not=|Q\cap W|$.

Hence, ${\cal M}$ is a resolving set. Finally, we obtain the bound by observing that
$|{\cal M}|=2e+1$. \qed \end{pf}

Let $\mathbb{F}_q$ be a finite field with $q$ elements, where $q$ is
a prime power. For a non-negative integer $n$, let $\mathbb{F}_q^{n}$
be an $n$-dimensional  vector space  over $\mathbb{F}_q$.
For a non-negative integer $m\leq n$, let $\left[[n]\atop m\right]$
be the set of all $m$-dimensional subspaces of $\mathbb{F}_q^{n}$. Then
the size of $\left[[n]\atop m\right]$ is
$\left[n\atop m\right]_q=\prod\limits_{i=n-m+1}^n(q^i-1)\big/\prod\limits_{i=1}^m(q^i-1).$

A {\it partition} of the vector space $V$ is a set ${\cal P}$ of
subspaces of $V$ such that any
non-zero vector is contained in exactly one element of ${\cal P}$.
If $T = \{\dim W \mid W\in{\cal P}\}$,
the partition $\cal P$ is said to be a $T$-{\it partition} of $V$.

Let $n=2e+1$ and $e\geq1$. The doubled Grassmann graph, denoted by $J_q(2e+1,e,e+1)$,
is a bipartite graph with bipartition $\left[[2e+1]\atop e\right]\cup\left[[2e+1]\atop e+1\right]$ such that
two vertices $P\in\left[[2e+1]\atop e\right]$ and $Q\in\left[[2e+1]\atop e+1\right]$
are adjacent if and only if $P\subseteq Q$.
Then $J_q(2e+1,e,e+1)$ is a distance-regular graph
with diameter $2e+1$ (see \cite{BCN}).
For any two vertices $P$ and $Q$, we have
$d(P,Q)=2i+|\dim P-\dim Q|$ if and only if $\dim(P\cap Q)=\min\{\dim P,\dim Q\}-i$.
For $e=1$, it is obvious that $\mu(J_q(3,1,2))= q(q+1)$.
For $e\geq2$, we obtain the following result.

\begin{thm}\label{thm1.2}
Let $e\geq2$. Then $\mu(J_q(2e+1,e,e+1))\leq(q^{2e+2}-1)/(q-1)$.
\end{thm}

\begin{pf}
By \cite[Lemma~3]{Beutel},  there exists an $\{e+1,e\}$-partition
${\cal P}=\{X,Y_1,\ldots,Y_{q^{e+1}}\}$ of $\mathbb{F}_q^{2e+1}$,
where $\dim X=e+1$ and $\dim Y_i=e,\,i=1,2,\ldots,q^{e+1}$. For a
fixed 1-dimensional subspace $U$ of $X$, let
$${\cal W}=\{W\subseteq U\oplus Y\mid \dim W=e, Y\in \{Y_1,\ldots,Y_{q^{e+1}}\}\},\quad
{\cal X}=\{W\subseteq X \mid \dim W=e\}.$$
To prove ${\cal M}={\cal W}\cup{\cal X}$ is a resolving set,
we only need to show that,
for any two distinct vertices $P$ and $Q$ in the same part,
there exists a vertex $W\in{\cal M}$ such that
$\dim(P\cap W)\not=\dim(Q\cap W)$.

Suppose $P,Q\in\left[[2e+1]\atop e+i\right]$, where $i=0$ or 1.

{\it Case} 1: $P\cap X\not=Q\cap X$. Suppose
$s+i=\dim (P\cap X)\leq\dim(Q\cap X)=r+i$. Then $1-i\leq r\leq e$.

{\it Case} 1.1: $s<r$. For $i=0$, let $W$ be an $e$-dimensional subspace of $X$
containing $Q\cap X$. Then $W\in{\cal M}$ and $$\dim(P\cap W)\leq s<r=\dim(Q\cap W).$$
For $i=1$, let  $W$ be an $e$-dimensional subspace of $X$ satisfying $\dim(P\cap W)=s$.
Then $W\in{\cal M}$ and $\dim(P\cap W)=s< e=\dim(Q\cap W).$

{\it Case} 1.2: $s=r$. Let $\{\alpha_1,\ldots,\alpha_r,\alpha_{r+i}\}$
be a basis for $Q\cap X$. Since $P\cap X\not=Q\cap X$,
there exists a $\beta\in (P\cap X)\backslash (Q\cap X)$ such that
$\{\alpha_1,\ldots,\alpha_r,\alpha_{r+i},\beta\}$ is linearly independent.
Extend this to a basis $\{\alpha_1,\ldots,\alpha_r,\alpha_{r+i},\beta,\gamma_1,\ldots,\gamma_{e-r-i}\}$
for $X$. Let $W$ be the $e$-dimensional subspace spanned by
$\{\alpha_1,\ldots,\alpha_r,\alpha_{r+i},\gamma_1,\ldots,\gamma_{e-r-i}\}$.
Then $W\in{\cal M}$. By $\beta\not\in W$ and $\beta\in P\cap X$,
$$\dim(P\cap W)<\dim(P\cap X)=\dim(Q\cap X)=r+i=\dim(Q\cap W).$$

{\it Case} 2: $P\cap X=Q\cap X$. Then there exists
a $Y\in \{Y_1,\ldots,Y_{q^{e+1}}\}$ such that $P\cap Y\not=Q\cap Y$.
Then $P\cap(U\oplus Y)\not=Q\cap(U\oplus Y)$.
Similar to the proof of Case~1, there exists a $W\in{\cal M}$
such that $\dim(P\cap W)\not=\dim(Q\cap W)$.

Hence, ${\cal M}$ is a resolving set. Finally, we obtain the bound by observing that
$$|{\cal M}|=\left[e+1\atop 1\right]_q+q^{e+1}\left[e+1\atop 1\right]_q=
\frac{q^{2e+2}-1}{q-1}.\eqno\Box
$$
\end{pf}

\section{Twisted Grassmann graphs}

For $e\geq2$,  let $H$ be a fixed $2e$-dimensional subspace of $\mathbb{F}_q^{2e+1}$,
\begin{eqnarray*}
{\cal B}_1&:=&\{W \; \hbox{a subspace of}\; \mathbb{F}_q^{2e+1}\mid \dim W=e+1, W\not\subseteq H\},\\
{\cal B}_2&:=&\{W \; \hbox{a subspace of}\; \mathbb{F}_q^{2e+1}\mid \dim W=e-1, W\subseteq H\}.
\end{eqnarray*}
The twisted Grassmann graph $\tilde{J}_q(2e+1,e)$ has the vertex set ${\cal B}_1\cup{\cal B}_2$,
and two vertices $P$ and
$Q$ are adjacent if and only if $\dim P+\dim Q-2\dim(P\cap Q)=2$.
Dam and Koolen \cite{Dam} constructed the graph, and showed that
$\tilde{J}_q(2e+1,e)$ is a distance-regular graph,
which is the first know family of non-vertex-transitive distance-regular graphs
with unbounded diameter.
For any two vertices $P$ and $Q$, $d(P,Q)=i$ if and only if $\dim P+\dim Q-2\dim(P\cap Q)=2i$.

\begin{thm}\label{thm1.1}
Let $e\geq2$. Then $\mu(\tilde{J}_q(2e+1,e))\leq (q^{2e}(q^{e+2}-q+1)-1)/(q-1)$.
\end{thm}

\begin{pf}
For $e\geq2$, by \cite[Lemma~3]{Beutel},  there exists an $\{e,1\}$-partition
${\cal P}=\{X_1,\ldots,X_{q^e+1},Y_1,\ldots,Y_{q^{2e}}\}$ of $\mathbb{F}_q^{2e+1}$,
where $\{X_1,\ldots,X_{q^e+1}\}$ is an $\{e\}$-partition of $H$
and  $\dim Y_i=1,\;i=1,\ldots,q^{2e}$.
For a fixed $(e+1)$-dimensional subspace $U$ of $H$, let
$${\cal W}=\{W\subseteq U\oplus Y\mid \dim W=e+1, W\not=U, Y\in \{Y_1,\ldots,Y_{q^{2e}}\}\},\;
{\cal X}=\{W\subseteq X \mid \dim W=e-1, X\in\{X_1,\ldots,X_{q^e+1}\}\}.$$
To prove
${\cal M}={\cal W}\cup{\cal X}$ is a resolving set, it is sufficient to
show that, for any two distinct vertices $P$ and $Q$,
there exists a vertex $W\in{\cal M}$ such that
$$
\left\{\begin{array}{ll}
\dim(P\cap W)\not=\dim(Q\cap W),&\hbox{if}\;P,Q\in{\cal B}_i;\\
\dim(P\cap W)-\dim(Q\cap W)\not=1,&\hbox{if}\;P\in{\cal B}_1,Q\in{\cal B}_2.
\end{array}\right.
$$

Firstly,  we assume that $P,Q\in{\cal B}_2$.
Similar to the proof of \cite[Proposition~7]{BM}, there exists
a $W\in{\cal X}$ such that $\dim(P\cap W)\not=\dim(Q\cap W)$.

\medskip
Secondly, we assume that $P,Q\in{\cal B}_1$.

{\it Case} 1: $P\cap H\not=Q\cap H$. Since $\{X_1,\ldots,X_{q^e+1}\}$
is a partition of $H$, there exists
an $X\in\{X_1,\ldots,X_{q^e+1}\}$ such that
$P\cap X\not= Q\cap X$. Suppose
$s=\dim (P\cap X)\leq\dim(Q\cap X)=r$. Then $0<r\leq e$.

{\it Case} 1.1: $r=e$. Then $X\subseteq Q$.
Since $\dim(Q\cap H)=e$, $Q\cap X=Q\cap H$.
Since $P\cap X\not= Q\cap X$,
$\dim(P\cap X)<e$. If $s=0$, take an $(e-1)$-dimensional $W$ of $X$. Then
$W\in{\cal M}$ and $0=\dim(P\cap W)<e-1=\dim(Q\cap W)$.
If $s\geq1$, let $W$ be an $(e-1)$-dimensional subspace of $X$ satisfying $\dim(P\cap W)=s-1$.
Then $W\in{\cal M}$ and
$$\dim(P\cap W)=s-1< e-1=\dim(Q\cap W).$$

{\it Case} 1.2: $s<r<e$. Let $W$ be an $(e-1)$-dimensional subspace of $X$
containing $Q\cap X$. Then $W\in{\cal M}$ and $$\dim(P\cap W)\leq s<r=\dim(Q\cap W).$$

{\it Case} 1.3: $s=r<e$. Let $\{\alpha_1,\ldots,\alpha_r\}$
be a basis for $Q\cap X$. Since $P\cap X\not=Q\cap X$,
there exists a $\beta\in (P\cap X)\backslash (Q\cap X)$ such that
$\{\alpha_1,\ldots,\alpha_r,\beta\}$ is linearly independent.
Extend this to a basis $\{\alpha_1,\ldots,\alpha_r,\beta,\gamma_1,\ldots,\gamma_{e-r-1}\}$
for $X$. Let $W$ be the $(e-1)$-dimensional subspace spanned by
$\{\alpha_1,\ldots,\alpha_r,\gamma_1,\ldots,\gamma_{e-r-1}\}$.
Then $W\in{\cal M}$. By $\beta\not\in W$ and $\beta\in P\cap X$,
$$\dim(P\cap W)<\dim(P\cap X)=\dim(Q\cap X)=r=\dim(Q\cap W).$$

{\it Case} 2: $P\cap H=Q\cap H$. Then there exists
a $Y\in \{Y_1,\ldots,Y_{q^{2e}}\}$ such that $P\cap Y\not=Q\cap Y$.
Then $P\cap(U\oplus Y)\not=Q\cap(U\oplus Y)$.
Suppose $s=\dim(P\cap(U\oplus Y))\leq\dim(Q\cap(U\oplus Y))=r$.
Then $2\leq s$ and $r\leq e+1$.

{\it Case} 2.1: $r=e+1$. Then $Q\subseteq U\oplus Y$ but $Q\not=U$,
which follows that $Q\in{\cal M}$ and $\dim(P\cap Q)<e+1$.

{\it Case} 2.2: $s<r<e+1$. Then exists an $(e+1)$-dimensional subspace $W$
of $U\oplus Y$ such that $\dim(P\cap W)=s-1<r-1\leq\dim(Q\cap W)$.
Since $P\cap U=Q\cap U$, $W\not= U$, which implies that $W\in{\cal M}$.

{\it Case} 2.3: $s=r<e+1$. Then $Q\cap (U\oplus Y)\not=Q\cap U$. Indeed, if $Q\cap (U\oplus Y)=Q\cap U$,
$P\cap U=Q\cap U$ by $U\subseteq H$.
Since $P\cap U\subseteq P\cap(U\oplus Y)$ and $\dim(P\cap U)\leq s\leq r=\dim(Q\cap U)$,
$P\cap U=P\cap(U\oplus Y)$, a contradiction. Let $\{\alpha_1,\ldots,\alpha_r\}$
be a basis for $Q\cap (U\oplus Y)$. Since $P\cap (U\oplus Y)\not=Q\cap (U\oplus Y)$,
there exists a $\beta\in (P\cap (U\oplus Y))\backslash (Q\cap (U\oplus Y))$ such that
$\{\alpha_1,\ldots,\alpha_r,\beta\}$ is linearly independent.
Extend this to a basis $\{\alpha_1,\ldots,\alpha_r,\beta,\gamma_1,\ldots,\gamma_{e+1-r}\}$
for $U\oplus Y$. Let $W$ be the $(e+1)$-dimensional subspace spanned by
$\{\alpha_1,\ldots,\alpha_r,\gamma_1,\ldots,\gamma_{e+1-r}\}$.
Then $W\not= U$ by $Q\cap (U\oplus Y)\not=Q\cap U$,
 which implies that $W\in{\cal M}$. By $\beta\not\in W$ and $\beta\in P\cap (U\oplus Y)$,
$$\dim(P\cap W)<\dim(P\cap (U\oplus Y))=\dim(Q\cap (U\oplus Y))=r=\dim(Q\cap W).$$

Next, we assume that $P\in{\cal B}_1,Q\in{\cal B}_2$.

{\it Case} 3.1: $Q\in{\cal M}$. Then $\dim(P\cap Q)-\dim Q\not=1$.

{\it Case} 3.2: $Q\not\in{\cal M}$. Since $\{X_1,\ldots,X_{q^e+1}\}$
is a partition of $H$ and $Q\subseteq H$, there exists
an $X\in\{X_1,\ldots,X_{q^e+1}\}$ such that
$e-2\geq s=\dim(Q\cap X)\geq1$. Let $\dim(P\cap X)=r$. Then $0\leq r\leq e$.

{\it Case} 3.2.1: $r>s+1$. If $r=e$, then $X\subseteq P$.
Let $\{\alpha_1,\ldots,\alpha_s\}$
be a basis for $Q\cap X$.
Extend this to a basis $\{\alpha_1,\ldots,\alpha_s,\gamma_1,\ldots,\gamma_{e-s}\}$
for $X$. Let $W$ be the $(e-1)$-dimensional subspace spanned by
$\{\alpha_2,\ldots,\alpha_s,\gamma_1,\ldots,\gamma_{e-s}\}$.
Then $W\in{\cal M}$. By $\alpha_1\not\in W$,
$\dim(Q\cap W)=s-1$. It follows that
$\dim(P\cap W)-\dim(Q\cap W)=e-s>1.$
Now let $r<e$. Suppose
$W$ is an $(e-1)$-dimensional
subspace of $X$ containing $P\cap X$. Then $W\in{\cal M}$ and $\dim(Q\cap W)\leq s$.
It follows that
$\dim(P\cap W)-\dim(Q\cap W)\geq r-s>1.$

{\it Case} 3.2.2: $r\leq s$. Let $W$ be an $(e-1)$-dimensional
subspace of $X$ containing $Q\cap X$. Then $W\in{\cal M},Q\cap X=Q\cap W$ and $\dim(P\cap W)\leq r$.
It follows that $\dim (P\cap W)-\dim(Q\cap W)\leq r-s<1$.

{\it Case} 3.2.3: $r= s+1$. If $Q\cap X\subseteq P\cap X$,
then let $\{\alpha_1,\ldots,\alpha_s,\alpha_{s+1}\}$
be a basis for $P\cap X$, where $\{\alpha_1,\ldots,\alpha_s\}$ is a basis for $Q\cap X$.
Extend this to a basis $\{\alpha_1,\ldots,\alpha_{s+1},\gamma_1,\ldots,\gamma_{e-s-1}\}$
for $X$. Let $W$ be the $(e-1)$-dimensional subspace spanned by
$\{\alpha_1,\ldots,\alpha_{s},\gamma_1,\ldots,\gamma_{e-s-1}\}$.
Then $W\in{\cal M}$. By $\alpha_{s+1}\not\in W$,
$\dim(P\cap W)=s$, which implies that
$\dim(P\cap W)-\dim(Q\cap W)=0.$
If $Q\cap X\not\subseteq P\cap X$, pick a basis
$\{\alpha_1,\alpha_2,\ldots,\alpha_{s}\}$ for $Q\cap X$. Then
there exists a $\beta\in (P\cap X)\backslash (Q\cap X)$ such that
$\{\alpha_1,\ldots,\alpha_{s},\beta\}$ is linearly independent.
Extend this to a basis $\{\alpha_1,\ldots,\alpha_{s},\beta,\gamma_1,\ldots,\gamma_{e-s-1}\}$
for $X$. Let $W$ be the $(e-1)$-dimensional subspace spanned by
$\{\alpha_1,\ldots,\alpha_s,\gamma_1,\ldots,\gamma_{e-s-1}\}$.
Then $W\in{\cal M}$. By $\beta\not\in W$ and $\beta\in P\cap X$,
$\beta\not\in P\cap W$. Hence $\dim(P\cap W)\leq s$ and
$\dim(P\cap W)-\dim(Q\cap W)\leq 0.$

\medskip
Hence, ${\cal M}$ is a resolving set. Finally, we obtain the bound by observing that
$$|{\cal M}|=(q^e+1)\left[e\atop 1\right]_q+q^{2e}\left(\left[e+2\atop
1\right]_q-1\right)=\frac{q^{2e}(q^{e+2}-q+1)-1}{q-1}.\eqno\Box
$$
\end{pf}

Babai \cite{Babai} obtained bounds on a parameter of primitive distance-regular graphs
which is equivalent to the metric dimension. A natural question is to compare our
result with those. For the case of the twisted Grassmann graph $\tilde{J}_q(2e+1,e)$,
Babai's most general bound (see \cite[Theorem~2.1]{Babai}) yields
$$
\mu(\tilde{J}_q(2e+1,e))<4\sqrt{\left[2e+1\atop e\right]_q}\log\left[2e+1\atop e\right]_q,
$$
while his stronger bound (see \cite[Theorem~2.4]{Babai}) yields
$$
\mu(\tilde{J}_q(2e+1,e))<2e\frac{\left[2e+1\atop e\right]_q}
{\left[2e+1\atop e\right]_q-M}\log\left[2e+1\atop e\right]_q,
$$
where
$$
M=\max_{1\leq j\leq e}q^{j^2}\left[e+1\atop
j\right]_q\left[e\atop j\right]_q.
$$
These bounds are difficult to evaluate exactly, so we conducted some
experiments using MATLAB to compare these bounds with the one
obtained in Theorem~\ref{thm1.2}.
Our experiments indicate that  our
constructive bound is an improvement on  Babai's general  bounds in
most of cases, but Babai's bound seems better than our bound for
small $q$.

\section*{Acknowledgement}
This research is partially supported by    NSF of China (10971052,
10871027),   NCET-08-0052, Langfang Teachers' College (LSZB201005),
  Hunan Provincial Natural Science Foundation
of China (09JJ3006), and   the Fundamental Research Funds for the
Central Universities of China.

\end{document}